\documentclass[12pt]{amsart}

\usepackage[latin2]{inputenc} 
\usepackage{amssymb,amsmath,latexsym,amscd,url}
\usepackage[mathscr]{eucal}
\usepackage[T1]{fontenc}

\theoremstyle{plain}
 
\newtheorem{thm}{Theorem}
\newtheorem{prop}{Proposition}

\input xy
\xyoption{all}

\newcommand{\slfrac}[2]{\left.#1\middle/_{#2}\right.}

\begin{document}
\title{Induced Cylindric Algebras of Choice Structures}
\author{Zolt\'an Moln\'ar}
\address{Department of Algebra, Budapest University of Technology and Economics, 1111 Egry J\'ozsef u. 1. Building H, $5^{\mathrm{th}}$ floor, Budapest, Hungary}
\email{mozow@math.bme.hu}
\urladdr{\url{http://www.math.bme.hu/~mozow/}}
\date{2010}

\begin{abstract}
One of the benefit properties implied by the extensionality axiom of Hilbert's epsilon calculus is that the calculus becomes complete with respect to the choice structures as semantics. Another implication of the axiom, discussed in the paper, is that an algebra is induced over the universe of the canonical model of a theory, which is isomorphic to a quotient algebra of the Lindenbaum--Tarski algebra of the theory. Especially, in the case of Boolean or monadic algebras, the canonical model of the theory of a sigma complete model is isomorphic to the algebra induced by the axiom of extensionality.
\end{abstract}
\maketitle
\section{Introduction}
The epsilon calculus was introduced by Hilbert and Bernays in \cite{Hilb}. Concerning the variable binding operator $(\varepsilon v_i)$, there are two main axiom-schemes. The axiom of transfinity (or the first epsilon axiom) implies the property that the term $(\varepsilon v_i)\varphi$ is a potential Henkin witness for the formula $\varphi$, i.e. the following formula is provable in every epsilon calculus
$$\varphi[t/v_i]\;\rightarrow\;\varphi[(\varepsilon v_i)\varphi/v_i]$$
where $[./v_i]$ is the operation of substitution, $t$ is any term and $v_i$ is a variable. The existential formula $(\exists v_i)\varphi$ and the universal formula $(\forall v_i)\varphi$ are defined by the substitutions $\varphi[(\varepsilon v_i)\varphi/v_i]$ and $\varphi[(\varepsilon v_i)\neg\varphi/v_i]$ respectively.\footnote{A brief introduction to the epsilon calculus can be found in \cite{Bour} pp. 36-44.} In Section 2, we shall prove that assuming the axiom of extensionality, the canonical model given by this method is unique, atomic and has an embedding property which is constructive in a sense. 

The extensionality axiom (or second epsilon axiom) is the following scheme
$$(\forall v_i)(\varphi\leftrightarrow \psi)\;\rightarrow\;(\varepsilon v_i)\varphi=(\varepsilon v_i)\psi$$
(where $\varphi,\psi$ are formulae and $v_i$ is a variable). On the one hand the extensionality axiom implies the model theoretic result of completeness.\footnote{See \cite{Leis}.} On the other hand there is a pure algebraic consequence of the axiom. The formula above induces a cylindric algebraic homomorphism from the algebra of formulae onto the epsilon terms. In Section 3, we shall introduce the algebra of epsilon terms, and we shall prove that it is a rich monadic algebra. In Section 4, we prove that in the case of Boolean or monadic algebras, the canonical model of the theory of a $\sigma$-complete model is isomorphic to the term algebra introduced in Section 3.
\section{Canonical models of extensional epsilon calculi}
In order to establish an appropriate environment for the algebraic approach, we deal with the canonical model\footnote{Cf. \cite{Hodg} p. 18.} of a complete and consistent theory in an epsilon language. Let $\mathsf{t}=(r_i, f_j, c_k)_{(i,j,k)\in I\times J\times K}$ be a similarity type and let $\mathscr{L}_\varepsilon$ be the freely generated language with respect to the operations $\neg,\vee, \varepsilon, \mathsf{t}$.\footnote{A detailed discussion of the language can be found in \cite{Monk} Ch. 29, p. 481} Let $\Gamma\subseteq \mathrm{Sent}(\mathscr{L}_\varepsilon)$ be a complete and consistent set of sentences. Let us define the canonical $\mathscr{L}_\varepsilon$-model $(\mathfrak{M},f)$ as follows.\footnote{The pair $(\mathfrak{M},f)$ will be a \emph{choice structure}, see \cite{Monk} p. 481.} Consider the set of epsilon terms of the single-variable formulae
$$\mathrm{Eps}_1=\{(\varepsilon v_i)\varphi\mid v_i\in\mathrm{Var}(\mathscr{L}_\varepsilon)\quad\mathrm{and}\quad \varphi\in\mathrm{Fm}_{v_i}(\mathscr{L}_\varepsilon)\}$$
and the equivalence relation $=_{\Gamma}$ over the set $\mathrm{Eps}_1$ as follows
$$t_1=_{\Gamma}t_2\quad\Longleftrightarrow\quad \Gamma\vdash t_1=t_2$$
for all $t_1,t_2\in\mathrm{Tm}(\mathscr{L}_\varepsilon)$.  Let the universe of the model $\mathfrak{M}$ be the set
$$M=\slfrac{\mathrm{Eps}_1}{\!=_{\Gamma}}$$
and let the interpretations of the relation, function and constant sings respectively be
\begin{align*}
(\slfrac{a_1}{=_\Gamma},...,\slfrac{a_l}{=_\Gamma})\in r_i^{\mathfrak{M}}&\quad\Leftrightarrow\quad \Gamma\vdash r_i(a_1/v_1,...,a_l/v_l)\\
f_j^{\mathfrak{M}}(\slfrac{a_1}{=_\Gamma},...,\slfrac{a_m}{=_\Gamma})&\quad=\quad \slfrac{(\varepsilon v_0)(v_0=f_j(a_1/v_1,...,a_m/v_m)}{=_{\Gamma}}\\
c_k^{\mathfrak{M}}&\quad=\quad \slfrac{(\varepsilon v_0)(v_0=c_k)}{=_{\Gamma}}
\end{align*}
The interpretations are well-defined, since they are independent of the choice of the representants. Then -- as it is well known from \cite{Leis} -- the axiom of extensionality and the axiom of transfinity allow us to define the \emph{canonical model}. Furthermore, this is a unique structure in the sense of the following proposition.

For the sake of simplicity we set 
$$\varphi(v_i)^{\Gamma}=\{(\slfrac{t}{=_{\Gamma}})\in M\mid \Gamma\vdash \varphi[t/v_i]\}$$
where $\varphi\in\mathrm{Fm}_{v_i}(\mathscr{L}_\varepsilon)$.
\begin{prop} Let $M$, $\mathfrak{M}$, $\Gamma$ be as above.
\begin{enumerate} 
\item[(a)] For all $\varphi, \psi\in\mathrm{Fm}_{v_i}(\mathscr{L}_\varepsilon)$ 
$$\mbox{if }\varphi(v_i)^{\Gamma}=\psi(v_i)^{\Gamma}\mbox{ then }\slfrac{((\varepsilon v_i)\varphi)}{=_{\Gamma}}\;=\;\slfrac{((\varepsilon v_i)\psi)}{=_{\Gamma}}$$
\item[(b)] There is a choice function $f$ such that for all $\varphi\in\mathrm{Fm}(\mathscr{L}_\varepsilon)$, $t\in\mathrm{Tm}(\mathscr{L}_\varepsilon)$ and valuation  $a=(\slfrac{a_1}{=_{\Gamma}},\slfrac{a_2}{=_{\Gamma}},...)$ of $\mathfrak{M}$
$$
(\mathfrak{M},f)\models\varphi[a]\quad\mbox{iff}\quad\Gamma\vdash\varphi[a_1/v_1,...,a_n/v_n]
$$
and
$$
\Gamma\vdash t^{\mathfrak{M}f}[a]=t[a_1/v_1,...,a_n/v_n]
$$
moreover, on the set 
$$\{ \{(\slfrac{t}{=_{\Gamma}})\in M\mid \Gamma\vdash \varphi[t/v_i]\}\mid \varphi\in\mathrm{Fm}_{v_i}(\mathscr{L}_\varepsilon)\}$$ 
$f$ unique.
\end{enumerate}
\end{prop}
\noindent\textsl{Proof.} (a) Let $s\overset{\circ}{=}(\varepsilon v_i)\neg(\varphi\leftrightarrow\psi)$. Since $s\in\mathrm{Eps}_1$ then by $\varphi(v_i)^{\Gamma}=\psi(v_i)^{\Gamma}$ we have $\Gamma\vdash\varphi[s/v_i]\leftrightarrow\psi[s/v_i]$. According to the definition of the universal quantifier
$$(\forall v_i)(\varphi\leftrightarrow\psi)\stackrel{\circ}{=}(\varphi\leftrightarrow\psi)[(\varepsilon v_i)\neg(\varphi\leftrightarrow\psi)/v_i]$$
therefore
$$\Gamma\vdash\varphi[s/v_i]\leftrightarrow\psi[s/v_i]\quad\mbox{iff}\quad\Gamma\vdash(\varphi\leftrightarrow\psi)[s/v_i] \quad\mbox{iff}\quad\Gamma\vdash(\forall v_i)(\varphi\leftrightarrow\psi)$$
The formula $(\forall v_i)(\varphi\leftrightarrow\psi)\rightarrow (\varepsilon v_i)\varphi=(\varepsilon v_i)\psi$ is an instance of the \emph{axiom of extensionality}
hence
$$\Gamma\vdash(\varepsilon v_i)\varphi=(\varepsilon v_i)\psi$$
holds and obviously $\slfrac{((\varepsilon v_i)\varphi)}{=_{\Gamma}}\;=\;\slfrac{((\varepsilon v_i)\psi)}{=_{\Gamma}}$.\\\\
(b) Let $f_*$ be a choice function such that $f_*(S)\in S$ if $S\in \mathrm{Sb}(M)\setminus\{ \emptyset\}$ and $f_*(\emptyset)\in M$. Let $f:\mathrm{Sb}(M)\to M$ be the function
$$f(S)=
\begin{cases}\slfrac{((\varepsilon v_i)\varphi)}{=_{\Gamma}} & \mbox{ if } S=\{(\slfrac{t}{=_{\Gamma}})\in M\mid \Gamma\vdash \varphi[t/v_i]\}\mbox{ for some }\varphi\in\mathrm{Fm}_{v_i}(\mathscr{L}_\varepsilon)\\
f_*(S) & \mbox{ otherwise }
\end{cases}$$
If $S\in \mathrm{Sb}(M)\setminus\{ \emptyset\}$ then $S=\{(\slfrac{t}{=_{\Gamma}})\in M\mid \Gamma\vdash \varphi[t/v_i]\}$  for some $\varphi\in\mathrm{Fm}_{v_i}(\mathscr{L}_\varepsilon)$. By the axiom of transfinity, we have $\vdash\varphi[t/v_i]\rightarrow \varphi[(\varepsilon v_i)\varphi/v_i]$, hence $\slfrac{(\varepsilon v_i)\varphi}{=_{\Gamma}}\in S$, $f$ is also a choice function and $(\mathfrak{M},f)$ is a choice structure in the language $\mathscr{L}_\varepsilon$. \\
The uniqueness can easily be shown by structural induction. $\square$
\begin{prop}
The canonical model of a complete and consistent theory is atomic.  
\end{prop}
\noindent\textsl{Proof.} Let $\Gamma\subseteq \mathrm{Sent}(\mathscr{L}_\varepsilon)$ be a complete and consistent set of sentences and let $(\mathfrak{M},f)$ be the canonical model of $\Gamma$. If $\slfrac{(\varepsilon v_i)\varphi}{=_{\Gamma}}\in M$ then $((\varepsilon v_i)\varphi)^{\mathfrak{M}f}=\slfrac{(\varepsilon v_i)\varphi}{=_{\Gamma}}$. Therefore $(\mathfrak{M},f)$ is atomic\footnote{Cf. \cite{Chan} Exercise 2.3.2., p. 107.}. $\square$ \\\\ 
Since the canonical model is atomic, it can be elementarily embedded into any model of its theory.\footnote{See \cite{Chan} Theorem 2.3.4., p. 99.} What is more, the epsilon terms give rise to the existence of an embedding which is defined straightforward by the values of the epsilon terms.\\\\  
We define the \emph{canonical injection} $\eta$ of the complete and consistent theory $\Gamma$. Let $(\mathfrak{N},g)$ be a model of $\Gamma$ and let us denote the canonical model of $\Gamma$ by $\mathfrak{Can}(\Gamma)$ and its universe by $\mathrm{Can}(\Gamma)$ then the canonical injection is 
$$\eta:\mathrm{Can}(\Gamma)\to N,\;\eta((\varepsilon v_i)\varphi)^{\mathfrak{Can}\,\Gamma}=((\varepsilon v_i)\varphi)^{\mathfrak{N}g}$$
Let us denote the canonical model of the theory $\mathrm{Th}(\mathfrak{N},g)$ by $\mathfrak{Can\,N}g$ and its universe by $\mathrm{Can}\,\mathfrak{N}g$.
\begin{prop}
If $(\mathfrak{N},g)$ is an $\mathscr{L}_\varepsilon$-model then the canonical injection $\eta:\mathrm{Can}\,\mathfrak{N}g\to N$ is an elementary embedding from  $\mathfrak{Can\,N}g$ to $(\mathfrak{N},g)$.
\end{prop}
\noindent\textsl{Proof.}  $\eta$ is a well-defined injection. Indeed, let us denote $\mathrm{Th}(\mathfrak{N},g)$ by $\Gamma$ and let $\slfrac{t}{=_{\Gamma}},\slfrac{s}{=_{\Gamma}}\in\mathrm{Can}\,\mathfrak{N}g$. By the definition of $\mathfrak{Can\,N}g$,
$$t^{\mathfrak{Can\,N}g}= s^{\mathfrak{Can\,N}g} \quad\mbox{ iff }\quad\mathrm{Th}(\mathfrak{N},g)\vdash t = s $$ 
holds, therefore if $\slfrac{t}{=_{\Gamma}}=\slfrac{s}{=_{\Gamma}}$, then $t^{\mathfrak{N}g}= s^{\mathfrak{N}g}$. Conversly, if $t^{\mathfrak{N}g}\ne s^{\mathfrak{N}g}$, then $(\mathfrak{N},g)\models t\ne s$
and $\mathrm{Th}(\mathfrak{N},g)\vdash t \ne s $.\\
$\eta$ is an elementary embedding. By definition, let $\varphi\in\mathrm{Fm}(\mathscr{L}_\varepsilon)$ and $a=(a_1/\!=_{\Gamma},a_2/\!=_{\Gamma},...)\in\mathrm{Val}(\mathfrak{Can\,N},g)$, then
\begin{align*}
\mathfrak{Can\,N}g\models\varphi[a]\quad&\mbox{iff}\quad\mathrm{Th}(\mathfrak{N},g)\vdash\varphi[a_1,...,a_n/x_1,...,x_n]\\
&\mbox{iff}\quad (\mathfrak{N},g)\models\varphi[a_1,...,a_n/x_1,...,x_n]\\
&\mbox{iff}\quad (\mathfrak{N},g)\models\varphi[\eta \circ a]\quad\square
\end{align*}
If $(\mathfrak{N}_1,g_1)$ and $(\mathfrak{N}_2,g_2)$ are $\mathscr{L}_\varepsilon$-models and $h:(\mathfrak{N}_1,g_1)\to(\mathfrak{N}_2,g_2)$ is a homomorhism such that 
$$\mathrm{Ran}(h\upharpoonright \mathrm{Eps}_1^{\mathfrak{N}_1g_1})\subseteq\mathrm{Eps}_1^{\mathfrak{N}_2g_2}$$
then let us define the function $h^{*}:\mathrm{Can}\,\mathfrak{N}_1g_1\to \mathrm{Can}\,\mathfrak{N}_2g_2$ by the relation
$$h^{*}(t^{\mathfrak{Can\,N}_1g_1})=s^{\mathfrak{Can\,N}_2g_2}\quad\mbox{whenever}\quad h(t^{\mathfrak{N}_1g_1})=s^{\mathfrak{N}_2g_2}$$
\begin{prop} Let $(\mathfrak{N}_1,g_1)$ and $(\mathfrak{N}_2,g_2)$ be $\mathscr{L}_\varepsilon$-models and let $h:(\mathfrak{N}_1,g_1)\to(\mathfrak{N}_2,g_2)$ be a homomorphism such that 
$$\mathrm{Ran}(h\upharpoonright \mathrm{Eps}_1^{\mathfrak{N}_1g_1})\subseteq\mathrm{Eps}_1^{\mathfrak{N}_2g_2}.$$ 
Then the canonical injection $\eta$ is a natural transformation in the sense that the following diagram commutes
$$\xymatrix{
\mathfrak{Can\,N}_1g_1\ar[r]^{h^*} \ar[d]_{\eta} &  \mathfrak{Can\,N}_2g_2 \ar[d]^{\eta} \\
(\mathfrak{N}_1,g_1)\ar[r]^{h}  & (\mathfrak{N}_2,g_2)}
$$
\end{prop}
\noindent\textsl{Proof.} By the definition of relation $h^*$ it follows that $h^*=\eta_1\circ h\circ \eta_2^{-1}$, where $\eta_i:\mathfrak{Can\,N}_ig_i\to (\mathfrak{N}_i,g_i)$ is the canonical injection. $\square$

\section{Canonical term algebras in extensional calculi}
The axiom of extensionality induces an algebra over the universe of the canonical model. Let  $(\mathfrak{M},f)$ be an $\mathscr{L}_{\varepsilon}$-model then $\mathbf{LT}\,\mathfrak{M}f$ denotes the Lindenbaum--Tarski algebra\footnote{Concerning the Lindenbaum--Tarski algebra and the cylindric set algebra of a first order model see \cite{Monk} p. 225.} of the formulae with the universe
$$\mathrm{LT}\,\mathfrak{M}f=\slfrac{\mathrm{Form}(\mathscr{L}_{\varepsilon})}{\mathfrak{M}\models.\leftrightarrow.}$$
The following proposition describes the connection between the neat-1-reduct $\mathbf{Nr}_1\,\mathbf{LT}\,\mathfrak{M}f$ of the cylindric algebra $\mathbf{LT}\,\mathfrak{M}f$ and the universe of the canonical model.
\begin{prop}\label{LT-Can}
If $(\mathfrak{M},f)\in\mathrm{Mod}(\mathscr{L}_{\varepsilon})$ then the function
$$\Phi:\mathrm{Nr}_1\,\mathbf{LT}\,\mathfrak{M}f\to \mathrm{Can}\,\mathfrak{M}f,\quad\Phi(\slfrac{\varphi}{\mathfrak{M}\models.\leftrightarrow.})=\slfrac{((\varepsilon v_0)\varphi)}{=_{\mathrm{Th}\,\mathfrak{M}f}}$$
is a surjection and if $\mathbf{Can}\,\mathfrak{M}f$ denotes the cylindric algebra generated by $\Phi$ then 
$$\frac{\mathbf{Nr}_1\,\mathbf{LT}\,\mathfrak{M}f}{\mathrm{Ker}\,\Phi}\cong\mathbf{Can}\,\mathfrak{M}f$$
\end{prop}
\noindent\textsl{Proof.} Let $\slfrac{((\varepsilon v)\varphi)}{=_{\mathrm{Th}\,\mathfrak{M}f}}\in \mathrm{Can}\,\mathfrak{M}f$ where $\varphi(v)\in\mathrm{Form}(\mathscr{L}_{\varepsilon})$. We can assume that $v_0$ is free for $v$ in $\varphi(v)$. Then  
$\slfrac{\varphi(v_0)}{\mathfrak{M}\models.\leftrightarrow.}\in\mathrm{Nr}_1\,\mathbf{LT}\,\mathfrak{M}f$ and $\mathfrak{M}\models(\varepsilon v)\varphi(v)=(\varepsilon v_0)\varphi(v_0)$. Furthermore, by the First Isomorphism Theorem there exists an isomorphism $\iota$ such that, the following diagram commutes
$$\xymatrix{
\mathbf{Nr}_1\,\mathbf{LT}\,\mathfrak{M}f\ar[r]^{\Phi}\ar[dr]_{\pi}  &  \mathbf{Can}\,\mathfrak{M}f \\
 & \frac{\mathrm{Nr}_1\,\mathbf{LT}\,\mathfrak{M}f}{\mathrm{Ker}\,\Phi} \ar[u]_\iota}
$$
where $\pi$ is the canonical projection. $\square$\\\\
We shall call $\mathbf{Can}\,\mathfrak{M}f$ the \emph{canonical term algebra} of $(\mathfrak{M},f)$. According to the following proposition, the canonical term algebra (as a monadic algebra) is rich\footnote{See \cite{Halm} p. 77.}. Moreover, the cylindrification $c_0$ of $\mathbf{Can}\,\mathfrak{M}f$ is a Boolean homomorphism from the Boolean reduct $\mathbf{Can}\,\mathfrak{M}f\upharpoonright\mathscr{L}_\varepsilon^\mathrm{BA}$ onto the trivial Boolean algebra $\mathbf{2}$.
\begin{prop} \label{rich} If $(\mathfrak{M},f)$ is a model, then 
\begin{enumerate} 
\item there exists a Boolean homomorphism $$b:\mathbf{Can}\,\mathfrak{M}f\upharpoonright\mathscr{L}_\varepsilon^\mathrm{BA}\to (\mathbf{Nr}_0\mathbf{Can}\,\mathfrak{M}f)\upharpoonright\mathscr{L}_\varepsilon^\mathrm{BA}$$ such that $b=c_0$ (where $c_0$ is the cylindrification)
\item  $\mathbf{Nr}_0\mathbf{Can}\,\mathfrak{M}f\upharpoonright\mathscr{L}_\varepsilon^\mathrm{BA} \cong\mathbf{2}$
\end{enumerate}
\end{prop}
\noindent\textsl{Proof.} (1) It immediately follows from the fact that every countable  monadic algebra is rich.\footnote{See \cite{Halm}, Part 3, Theorem 2, p. 77.} (2) $\mathbf{Nr}_0\,\mathbf{LT}\,\mathfrak{B}f\upharpoonright\mathscr{L}_\varepsilon^\mathrm{BA}$ is the 0-Lindenbaum--Tarski algebra of a complete and consistent theory, therefore it is the two-element Boolean algebra $\mathbf{2}$. Then claim (2) follows from the previous fact and Proposition \ref{LT-Can}.
$\square$
\section{Canonical term algebras of $\mathsf{BA}$s and $\mathsf{CA_1}$s}
When the model $(\mathfrak{M},f)$ is a $\sigma$-complete Boolean algebra, there is a close algebraic relationship between the cylindric set algebra of $(\mathfrak{M},f)$ and the canonical model $\mathfrak{Can}\,\mathfrak{M}f$. The theorem below describes the main connection.
\begin{thm} \label{main} Let $(\mathfrak{B},g)\in\mathsf{BA}$ and $(\mathfrak{M},f)\in\mathsf{CA}_1$ such that  $\mathfrak{Can}\,\mathfrak{B}g$ and $\mathfrak{Can}\,\mathfrak{M}f$ are finite algebras. Then 
\begin{enumerate}
\item  
$\mathfrak{Can}\,\mathfrak{B}g\cong\mathbf{Can}\,\mathfrak{B}g\upharpoonright\mathscr{L}_\varepsilon^\mathrm{BA}$
\item $\mathfrak{Can}\,\mathfrak{M}f\cong\mathbf{Can}\,\mathfrak{M}f\quad \mbox{iff} \quad \mathbf{Nr}_0\,\mathfrak{Can}\,\mathfrak{M}f\cong\mathbf{2}$
\end{enumerate}
\end{thm}
\noindent\textsl{Proof.} (1) By definition, universes of $\mathfrak{Can}\,\mathfrak{B}f$ and $\mathbf{Can}\,\mathfrak{B}f\upharpoonright\mathscr{L}_\varepsilon^\mathrm{BA}$ are the same. It is known that all finite Boolean algebras with the same cardinality are isomorphic.\footnote{Cf. \cite{Monk} Corollary 9.32., p. 152.} Hence, $$\mathfrak{Can}\,\mathfrak{M}f\upharpoonright\mathscr{L}_\varepsilon^\mathrm{BA}\cong\mathbf{Can}\,\mathfrak{M}f\upharpoonright\mathscr{L}_\varepsilon^\mathrm{BA}$$
(2) According to (1) the Boolean reducts of $\mathfrak{Can}\,\mathfrak{M}f$ and $\mathbf{Can}\,\mathfrak{M}f$ are isomorphic. By Proposition \ref{rich}, the cylindrification $c_0$ is a Boolean homomorphism. Hence, $\mathbf{Nr}_0\,\mathfrak{Can}\,\mathfrak{M}f\cong \mathbf{Nr}_0\,\mathbf{Can}\,\mathfrak{M}f$ iff $\mathbf{Nr}_0\,\mathfrak{Can}\,\mathfrak{M}f\cong\mathbf{2}$.
$\square$\\\\
As an application, we prove two sufficient conditions concerning the existence of the algebraic connection between the cylindric set algebra of the model and the canonical model.
\begin{prop}  Let $(\mathfrak{B},g)\in\mathsf{BA}$ and $(\mathfrak{M},f)\in\mathsf{CA}_1$.
\begin{enumerate}
\item  If $(\mathfrak{B},g)$ is $\sigma$-complete, then $$\frac{\mathbf{Cs}\,(\mathfrak{B},g)\upharpoonright\mathscr{L}_\varepsilon^\mathrm{BA}}{\mathrm{Ker}\,\Phi}\cong\mathfrak{Can}\,\mathfrak{B}g$$
\item  If $(\mathfrak{M},f) \upharpoonright\mathscr{L}_\varepsilon^\mathrm{BA}$ is $\sigma$-complete and $\mathbf{Nr}_0(\mathfrak{M},f)\cong\mathbf{2}$, then $$\frac{\mathbf{Cs}\,(\mathfrak{M},f)}{\mathrm{Ker}\,\Phi}\cong\mathfrak{Can}\,\mathfrak{M}f$$
\end{enumerate}
where $\Phi:\mathrm{Cs}\,(\mathfrak{M},f)\to \mathrm{Can}\,\mathfrak{M}f,\quad\Phi(\varphi^{\mathfrak{M}f})=\slfrac{((\varepsilon v_0)\varphi)}{=_{\mathrm{Th}\,\mathfrak{M}f}}$.\footnote{Cf. Ibid. Proposition \ref{LT-Can}.}
\end{prop}
\noindent\textsl{Proof.} It is clear that $\mathbf{Cs}\,(\mathfrak{M},f)\cong\mathbf{LT}\,\mathfrak{M}f$. The $\sigma$-completeness implies that, the canonical model is isomorphic to a subalgebra of $(\mathfrak{M},f)$ which is a power set algebra\footnote{See \cite{Kopp}, p. 31, CH. 1, Corollary 2.7.}. Hence, the canonical model is finite and by Proposition \ref{LT-Can} and Theorem \ref{main} it follows that the mentioned quotient algebra is isomorphic to the canonical model. $\square$

\end{document}